\documentclass{llncs}
\usepackage{hyperref}
\usepackage{epsfig,graphicx,amsfonts,amssymb}

\spnewtheorem{observation}{Observation}{\bfseries}{\itshape}
\spnewtheorem{clm}{Claim}{\bfseries}{\itshape}
\def\bsquareforqed{\rule{0.6em}{0.6em}}
\def\bqed{\ifmmode\bsquareforqed\else{\unskip\nobreak\hfil
\penalty50\hskip1em\null\nobreak\hfil\bsquareforqed
\parfillskip=0pt\finalhyphendemerits=0\endgraf}\fi}
\newcommand{\ds}{\displaystyle}  
\newcommand{\comp}[1]{\overline{#1}}  

\newcommand{\numed}[1]{|| #1 ||}
\newcommand{\cub}{\mbox{\textnormal{cub}}}
\newcommand{\sph}{\mbox{\textnormal{sph}}}
\newcommand{\boxi}{\mbox{\textnormal{box}}}
\newcommand{\lt}[1]{\ensuremath{#1^{-}}}
\newcommand{\rt}[1]{\ensuremath{#1^{+}}}
\newcommand{\ceil}[1]{\lceil #1 \rceil}
\newcommand{\ignore}[1]{}
\newcommand{\floor}[1]{\lfloor #1 \rfloor}
\newcommand{\vb}{N} 
\newcommand{\sn}{N_S} 

\newcommand{\igr}{\mathcal{I}}
\newcommand{\nscomm}[1]{\textbf{#1}}
\newcommand{\imin}{I_{min}}

\begin{document}
\title{Lower Bounds for Boxicity}

\author{Abhijin Adiga\inst{1}, L. Sunil Chandran\inst{2}, and Naveen
Sivadasan\inst{3}}
\institute{
Network Dynamics and Simulation Science Laboratory,
Virginia Bioinformatics Institute, Virginia Tech,
Blacksburg, VA 24061, USA.\\
email: abhijin@vbi.vt.edu
\and   
Department of Computer Science and Automation, Indian Institute
of Science, Bangalore -- 560012, India. \\email: sunil@csa.iisc.ernet.in
\and
Department of Computer Science and Engineering, Indian Institute of
Technology, Hyderabad -- 502205, India. \\email: nsivadasan@iith.ac.in}
\date{}
\maketitle
\begin{abstract}
An axis-parallel $b$-dimensional box is a Cartesian product $R_1\times
R_2\times\cdots\times R_b$ where $R_i$ is a closed interval of the
form $[a_i,b_i]$ on the real line. For a graph $G$, its \emph{boxicity}
$\boxi(G)$ is the minimum dimension $b$, such that $G$ is representable
as the intersection graph of boxes in $b$-dimensional space.  Although
boxicity was introduced in 1969 and studied extensively, there are no
significant results on lower bounds for boxicity. In this paper, we
develop two general methods for deriving lower bounds. Applying these
methods we give several results, some of which are listed below:
\begin{enumerate}
\item The boxicity of a graph on $n$ vertices with no universal vertices
and minimum degree $\delta$ is at least $n/2(n-\delta-1)$.
\item Consider the $\mathcal{G}(n,p)$ model of random graphs. Let 
$ p\le1-  \frac{ 40 \log n}{n^2}$. Then with high probability,
$\boxi(G)=\Omega(np(1-p))$. On setting $p=1/2$ we immediately infer
that almost all graphs have boxicity $\Omega(n)$. Another consequence of
this result is as follows: 
Let $m$ be an integer such that $m\le c n^2/3$, where $c < 1$ is any positive constant.
Then for almost all graphs $G$ on $n$ vertices and $m$ edges, $\boxi(G) = \Omega(m/n)$.
\item Let $G$ be a connected $k$-regular graph on $n$ vertices. Let
$\lambda$ be the second largest eigenvalue in absolute value of the
adjacency matrix of $G$. Then, the boxicity of $G$ is at least
\(
\left(\frac{k^2/\lambda^2}{\log\left(1+k^2/\lambda^2\right)}\right)\left(\frac{n-k-1}{2n}\right).
\)
\item The boxicity of random $k$-regular graphs on $n$ vertices (where $k$ is
fixed) is $\Omega(k/\log k)$.
\item
Consider balanced bipartite graphs with $2n$ vertices and $m$ edges such
that $m\le c n^2/3$, where $c < 1$ is any positive constant.  Then for
almost all balanced bipartite graphs $G$ on $2n$ vertices and $m$ edges,
$\boxi(G) = \Omega(m/n)$.
\end{enumerate}
\end{abstract}
\section{Introduction}
Let $G(V,E)$ be a simple undirected graph with vertex set $V$ and
edge set $E$. Let $\mathcal{F} = \{S_x\subseteq U : x\in V\}$ be a family of
subsets of a universe $U$, where $V$ is an
index set. The intersection graph $\Omega(\mathcal{F})$ of $\mathcal {F}$ has
 $V$ as
vertex set, and two distinct vertices $x$ and $y$ are adjacent if and only
if $S_x\cap S_y\ne\emptyset$.
An axis-parallel $b$-dimensional box is a Cartesian
product $R_1\times R_2\times\cdots\times R_b$ where $R_i$ is a closed
interval of the form $[a_i,b_i]$ on the real line. For a graph $G$,
its \emph{boxicity} $\boxi(G)$ is the minimum dimension $b$, such
that $G$ is representable as the intersection graph of boxes in
$b$-dimensional space.  Roberts \cite{recentProgressesInCombRoberts}
introduced the concept of boxicity. Its applications include niche overlap
(competition) in ecology and problems of fleet maintenance in operations
research. Cozzens \cite{phdThesisCozzens} showed that computing the
boxicity of a graph is NP-hard. This was later strengthened by Yannakakis
\cite{complexityPartialOrderDimnYannakakis} and finally, by Kratochvil
\cite{specialPlanarSatisfiabilityProbNPKratochvil}, who showed that
deciding whether boxicity of a graph is at most 2 itself is NP-complete.

The class of graphs with boxicity at most 1 is exactly the popular class
of \emph{interval graphs}. A graph $G(V,E)$ is an interval graph if
there exists a mapping, say $f$, from $V$ to intervals on the real line
such that two vertices in $G$ are adjacent if and only if the intervals
to which they are mapped overlap. $f$ is called an \emph{interval
representation} of $G$.
Given a graph $G(V,E)$, let $\igr=\{I_i(V,E_i),1\le i\le k\}$ be a set of
$k$ interval graphs such that $E=E_1\cap E_2\cap\cdots\cap E_k$. Then, we
say that $\mathcal{I}$ is an \emph{interval graph representation} of $G$.
The concept of boxicity can be equivalently formulated using the interval
graph representation as follows: 
\begin{theorem}
\textnormal{\textbf{(Roberts \cite{recentProgressesInCombRoberts})}}
The minimum $k$ such that there exists an interval graph representation
of $G$ using $k$ interval graphs  $\igr=\{I_i(V,E_i),1\le i\le k\}$
is the same as $\boxi(G)$.
\end{theorem}

\subsection{Some Known Bounds on Boxicity}
There have been several attempts to establish upper bounds on
boxicity, especially for graphs with special structures. Roberts,
in his seminal work \cite{recentProgressesInCombRoberts}, proved that
the boxicity of a complete $k$-partite graph is $k$, thereby showing
that boxicity of a graph with $n$ vertices cannot exceed $\floor{n/2}$.
Chandran and Sivadasan \cite{boxicityTreewidthSunilNaveen} showed that
$\boxi(G)\le\mbox{tw}(G)+2$, where, $\mbox{tw}(G)$ is the treewidth of
$G$. Chandran et al \cite{boxicityMaxDegreeSunilNaveenFrancis} proved
that $\boxi(G)\le\chi(G^2)$ where, $G^2$ is the supergraph of $G$ in
which two vertices are adjacent if and only if they are at a distance
of at most 2 in $G$, and $\chi(G^2)$ is the chromatic number of $G^2$. From
this result, they inferred that $\boxi(G)\le2\Delta^2$, where $\Delta$
is the maximum degree of $G$. Scheinerman \cite{phdThesisScheinerman} showed
that the boxicity of outer planar graphs is at most 2. Thomassen
\cite{intervalRepPlanarGraphsThomassen} proved that the boxicity of
planar graphs is at most 3. In \cite{computingBoxicityCozzensRoberts},
Cozzens and Roberts studied the boxicity of split graphs.

In contrast, the literature provides few results regarding lower
bounds on boxicity. Even ad hoc constructions that achieve high
boxicity are rare.  In \cite{recentProgressesInCombRoberts}, to prove
that a complete $k$-partite graph has boxicity at least $k$, Roberts
uses the pigeon-hole principle in conjunction with the fact that an
interval graph does not contain an induced $C_4$. As a consequence,
he shows that the boxicity of a $\ceil{n/2}$-partite graph is at least
$\floor{n/2}$. In \cite{forbiddenTrotter} Trotter has characterized
all graphs with boxicity $\floor{n/2}$. Similar arguments are used by
Cozzens and Roberts \cite{computingBoxicityCozzensRoberts} to show that
the boxicity of the complement of a cycle or a path of length $n$, is
at least $\ceil{n/3}$. Motivated by the concept of boxicity, McKee and
Scheinerman \cite{chordalityMcKeeScheinerman} study a parameter called
\emph{chordality} of a graph $G$, $\mbox{chord}(G)$ (a better name
would be chordal dimension). It is the minimum $k$ such that $G$ can be
expressed as the intersection of $k$ chordal graphs. Since, every interval
graph is chordal, it follows that $\boxi(G)\ge\mbox{chord}(G)$. To
obtain lower bounds for chordality, the authors use the property that
if a graph has chordality at most $k$, then $G$ contains a vertex whose
neighbours induce a subgraph of chordality at most $k-1$. Using this,
they have managed to show that the chordality and hence the boxicity
of a bipyramid graph is at least $3$. Unfortunately, it does not
look like this method is adequate to give strong and general lower bounds. In
\cite{boxicityTreewidthSunilNaveen}, Chandran and Sivadasan provide a
specialized construction to show that for any integer $k\ge1$, there
are graphs with treewidth at most $t+\sqrt{t}$ whose boxicity is at
least $t-\sqrt{t}$.

It is interesting to note that coloring problems on low boxicity
graphs were considered as early as 1948 \cite{problem56Bielecki}.
Kostochka \cite{coloringIntersectionGraphsKostochka} provides an
extensive survey on colouring problems of intersection graphs. In
\cite{findingImaiAsano,optimalPackingFowlerPatersonTanimoto}, the
complexity of finding the maximum independent set in bounded boxicity
graphs is considered.
Researchers have also tried to generalize or extend the concept of boxicity
in various ways. The poset boxicity \cite{posetBoxicityTrotterWest}, the
rectangle number \cite{rectNumHypercubesChangeWest}, grid dimension
\cite{gridIntBoxicityBellantoniEtal}, circular dimension
\cite{circDimensionFeinberg,noteCircDimensionShearer}, and the boxicity of
digraphs \cite{intNumBoxDigraphs} are some examples.

\subsection{Our Results}
In this paper, we present two methods to obtain lower bounds on
boxicity. The underlying idea of both the methods is to make use of the
vertex isoperimetric properties of the given graph (sometimes 
modified to suit our purpose).

Most of the already known lower
bounds can be re-derived using these methods.  Applying these methods
we derive several new results, some of which are listed below:
\begin{enumerate}
\item The boxicity of a graph on $n$ vertices with no universal vertices and
minimum degree $\delta$ is at least $n/2(n-\delta-1)$ (Theorem
\ref{thm:universal}, Section \ref{sec:consMainThm}). 
\begin{remark}
There are two
parameters, namely \emph{cubicity} and \emph{sphericity} which are closely
related to boxicity. There have been several instances in the literature when these
two parameters were compared with each other.  Cubicity (Sphericity)
of $G$ denoted as $\cub(G)$ ($\sph(G)$) is the minimum dimension $b$, such
that $G$ can be represented as the intersection graph of axis-parallel
unit cubes (unit spheres) in $b$-dimensional space. Since cubicity is
a stricter notion of boxicity, $\cub(G)\ge\boxi(G)$. Apparently, there
is no such relationship between cubicity (or boxicity) and sphericity.

Havel \cite{phdThesisHavel} showed that there
are graphs with sphericity 2 and arbitrarily high cubicity while Fishburn
\cite{onTheSphCubGraphsFishburn} constructed some graphs of cubicity at
most 3 with sphericity greater than their cubicity. Maehara et al
\cite{embeddingTreesMaeharaEtal} proved that the sphericity of the
complement of a tree is at most 3. Using our result mentioned above, we can
easily infer that if $G$ is the complement of a bounded degree tree, then
$\boxi(G)=\Omega(n)$. Thus we have a large number of graphs that have
arbitrarily higher value of boxicity (not just cubicity!) than their
sphericity. This is a much stronger result compared to that of Havel's, who
refers to the class of star graphs (Consider the five--pointed star
graph $K_{1,5}$: $\cub(K_{1,5})=3$ while $\sph(K_{1,5})=2$. But the
boxicity of any star graph is 1).
\end{remark}
\item Consider the $\mathcal{G}(n,p)$ model of random graphs. Let
$ p\le1-\frac{40 \log n}{n^2}$.
 Then with high probability,  
$\boxi(G)=\Omega(np(1-p))$ (Theorem \ref{thm:randThm}, Section
\ref{sec:rand1}).

On setting $p=1/2$ we immediately infer that almost all graphs have
boxicity $\Omega(n)$. Another consequence of this result is as follows:
Let $m$ be an integer such that $m\le c n^2/3$, where $c < 1$ is any positive constant.
Then for almost all graphs $G$ on $n$ vertices and $m$ edges, $\boxi(G) = \Omega(m/n)$.
\item Let $G$ be a connected $k$-regular graph on $n$ vertices. Let
$\lambda$ be the second largest eigenvalue in absolute value of the
adjacency matrix of $G$. Then, the boxicity of $G$ is at least
\(
\left(\frac{k^2/\lambda^2}{\log\left(1+k^2/\lambda^2\right)}\right)\left(\frac{n-k-1}{2n}\right)
\)
(Theorem \ref{thm:spec}, Section \ref{sec:spec}).

Bilu and Linial \cite{sphericityBiluLinial} have a similar result for
sphericity of regular graphs. They prove that an $\epsilon n$-regular
graph of order $n$, with bounded diameter has sphericity
$\Omega(n/(\lambda+1))$, where $0\le\epsilon\le\frac{1}{2}$ is a constant.
Note that if $G$ is an $\epsilon n$-regular graph, the lower bound for
boxicity that we get is almost the square of the Bilu-Linial bound for
sphericity.
\item The boxicity of random $k$-regular graphs on $n$ vertices (where $k$ is
fixed) is $\Omega(k/\log k)$ (Corollary \ref{cor:specAlmost}, Section
\ref{sec:spec}).


\item
Consider balanced bipartite graphs with $2n$ vertices and $m$ edges such
that $m\le c n^2/3$, where $c < 1$ is any positive constant.  Then for
almost all balanced bipartite graphs $G$ on $2n$ vertices and $m$ edges,
$\boxi(G) = \Omega(m/n)$
(Corollary \ref{cor:randBipartiteThmEdges}, Section \ref{sec:randbipartite}). 

\end{enumerate}

\section{Preliminaries}\label{sec:prelims}
For a graph $G(V,E)$, $\numed{G}$ denotes the number of edges in $G$.
The complement of $G$ is denoted by $\comp{G}$. A balanced
($A,B$)-bipartite graph is a bipartite graph with $V=A\uplus B$
and $|A|=|B|$. $\delta(G)$ and $\Delta(G)$ are the minimum degree and
maximum degree of $G$ respectively.

\subsection{Interval Graphs}\label{sec:intProp}
Suppose $I(V,E)$ is an interval graph and $f$ is an
interval representation of $I$. For a vertex $u$, let $l_f(u)$
and $r_f(u)$ denote the real numbers corresponding to the left end-point
and right end-point respectively of the interval $f(u)$. When there is
no ambiguity regarding the interval representation under consideration,
we shall discard the subscript $f$ and use the abbreviated forms $l(u)$
and $r(u)$ respectively. Further, we refer to $l(u)$ as the ``left
end-point" and $r(u)$ as the ``right end-point".

We can assume without loss of generality that for any interval graph $G$,
there is an interval representation such that all the interval end points
map to distinct points on the real line.

\paragraph{\textnormal{\textbf{Induced vertex numbering:}}} Given an
interval representation of $I$ with distinct end-points, we define the
induced vertex numbering $\eta(\cdot)$ as a numbering of the vertices in the
increasing order of their right end-points, i.e. for any two distinct
vertices $u$ and $v$, $\eta(u)<\eta(v)$ \nscomm{if and only if} $r(u)<r(v)$.

\begin{definition}
A \textbf{minimum interval supergraph} of a graph $G(V,E)$, is an interval
supergraph of $G$ with vertex set $V$ and with least number of edges
among all interval supergraphs of $G$.
\end{definition}

\subsection{Neighbourhoods}\label{sec:neighbours}
Given a subset of vertices, say $X$, we define several vertex
neighbourhoods of $X$.
\begin{enumerate}
\item{\textbf{Vertex-boundary}} $\vb(X,G) = \{u\in V-X|\exists v\in
X\mbox{~with~}uv\in E\}$. The term vertex-boundary is
borrowed from \cite{globalMethodsHarper}.  
\item{\textbf{Strong vertex-boundary}} $\sn(X,G)=\left\{u\in V-X| uv\in E,
\forall v\in X\right\}$.
\item{\textbf{Vertex neighbourhood}} $N'(X,G) = \{u|\exists v\in
X\mbox{~with~}uv\in E\}$.
\item{\textbf{Closed vertex neighbourhood}} $N[X,G] = X \cup N(X,G)$.
\end{enumerate}
When $X$ is a singleton set, say $X=\{v\}$, we will use the notations
$\vb(v,G)$, $\sn(v,G)$, $N'(v,G)$ and $N[v,G]$ respectively.
When there is no ambiguity regarding the graph under consideration,
we discard the second argument and simply denote the neighbourhoods
as $\vb(X)$, $\sn(X)$, $N'(X)$ and $N[X]$ respectively.

\subsubsection{Probability Spaces of Random Graphs}
\begin{enumerate}
\item $\mathcal{G}(n,p)$: The probability space of graphs on $n$
vertices in which the edges are chosen independently with probability
$p$.  
\item $\mathcal{G}(n,m)$: The probability space of graphs on $n$
vertices having $m$ edges, in which the graphs have the same probability.
\item $\mathcal{G}_R(n,k)$: The probability space of all simple $k$-regular
graphs on $n$ vertices with each graph having the same probability.
\item $\mathcal{G}_B(2n,p)$: The probability space of all balanced
bipartite graphs on $2n$ vertices in which the edges are chosen
independently with probability $p$.  
\item $\mathcal{G}_B(2n,m)$: The probability space of all balanced
bipartite graphs on $2n$ vertices and $m$ edges, in which the graphs
have the same probability.
\end{enumerate}
If $X$ is a random variable, $E(X)$ denotes the expectation of $X$. 

\subsubsection{Standard Asymptotic Notations:}
(From \cite{probabilisticMethodAlonSpencerErdos}) For two functions
$f$ and $g$, we write $f=\mathcal{O}(g)$ if $f\le c_1g+c_2$ for all
possible values of the variables of the two functions, where $c_1,c_2$
are absolute constants. We write $f=\Omega(g)$ if $g=\mathcal{O}(f)$
and $f=\Theta(g)$ if $f=\mathcal{O}(g)$ and $f=\Omega(g)$. If the limit
of the ratio $f/g$ tends to $0$ as the variables of the functions tend
to infinity, we write $f=o(g)$. We use $f=\omega(g)$ if $g=o(f)$.

\section{The First Method: Based on Minimum Interval Supergraph}
\label{sec:firstMethod}
Given any positive integer $k$, the \emph{vertex-isoperimetric problem}
\cite{globalMethodsHarper} is to minimize $|\vb(X,G)|$ over all
$X\subseteq V$ such that $|X|=k$. Let 
\begin{equation}\label{eqn:bvi}
b_v(k,G)=\min_{\stackrel{X\subseteq V}{|X|=k}}|\vb(X,G)|.
\end{equation}
A brief introduction to isoperimetric problems can be found in
\cite{combinatoricsBollabas}. Harper \cite{globalMethodsHarper} gives
a detailed treatment of the vertex-isoperimetric problem.

Let $c_v(k,G)=\max_{\stackrel{X\subseteq V}{|X|=k}}|\sn(X,G)|$.
It is easy to see that $\sn(X,\comp{G})=V-X-\vb(X,G)$. From this we infer
the following: 
\begin{eqnarray}
c_v(k,\comp{G})=\max_{\stackrel{X\subseteq
V}{|X|=k}}|V-X-\vb(X,G)|&&=n-k-\min_{\stackrel{X\subseteq
V}{|X|=k}}|\vb(X,G)|\nonumber\\
&&=n-k-b_v(k,G)\label{eqn:cvk}.
\end{eqnarray}
\begin{observation}\label{obs:cvObs1}
Let $i$ and $j$ be two positive integers such that $i>j$. Then,
$c_v(i,G)\le c_v(j,G)$.
\end{observation}
\begin{proof}
If $c_v(i,G)=k$, it implies that there exists a pair of sets
$X,Y\subseteq V$ with $|X|=i$ and $|Y|=k$ such that $Y$ is the strong
vertex boundary of $X$. For any set $A\subseteq X$ such that $|A|=j$,
$Y\subseteq\sn(A,G)$. Therefore, $c_v(j,G)\ge k$.
\end{proof}
\begin{observation}\label{obs:cvObs2}
Suppose $i$ and $j$ are non-negative integers. If $c_v(i,G)=j$, then,
for any integer $k>j$, $c_v(k,G)<i$.
\end{observation}
\begin{proof}
If $c_v(k,G)\ge i$, then, it implies that there exists a pair of sets
$X,Y\subseteq V$ with $|X|=k$ and $|Y|\ge i$, such that $Y=\sn(X,G)$. Then
clearly $X\subseteq \sn(Y,G)$. It follows that $c_v(i,G)\ge c_v(|Y|,G)\ge
k>j$, contradicting the assumption that $c_v(i,G)=j$.
\end{proof}

\subsection{Lower Bound for Boxicity}
\begin{lemma}\label{lem:iminLem} 
Let $G(V,E)$ be any non-complete graph. Let $\imin$ be a minimum interval
supergraph of $G$. Then,
\begin{equation}\label{eqn:gi}
\boxi(G)\ge \numed{\comp{G}}/\numed{\comp{\imin}}.
\end{equation}
\end{lemma}
\begin{proof}
Let $\igr=\{I_1,I_2,\ldots,I_k\}$ be an interval graph representation
of $G$. Since each graph $I\in \igr$ is an interval supergraph of $G$,
$\numed{I}\ge \numed{\imin}$, $\forall I\in \igr$. Hence, in any $I\in
\igr$, at most $\numed{\comp{\imin}}$ edges can be absent. In $G$,
$\numed{\comp{G}}$ edges are absent. Since, an edge absent in $G$
should be absent in at least one interval graph in $\igr$, there should
be at least $\numed{\comp{G}}/\numed{\comp{\imin}}$ interval graphs in
$\igr$. Hence proved.\bqed
\end{proof}

\begin{theorem}\label{thm:cvithm}
Let $G$ be a non-complete graph with $n$ vertices. Then,
\begin{equation}\label{eqn:cvi}
\boxi(G)\ge \frac{\numed{\comp{G}}}{\ds\sum_{i=1}^{n-1}c_v(i,\comp{G})}.
\end{equation}
\end{theorem}
\begin{proof}
In view of Lemma \ref{lem:iminLem}, it is enough to show that
$\numed{\comp{\imin}}\le\sum_{i=1}^{n-1}c_v(i,\comp{G})$. Consider
an interval representation with distinct end-points for $\imin$
and, let $\eta$ be the induced vertex numbering (see Section
\ref{sec:intProp} for definition). Given an integer $k$, $1\le k\le n$, let
\[
S_k(\eta)=\eta^{-1}\left(\{1,2,\ldots,k\}\right)=\left\{v\in V|\eta(v)\le k\right\},
\]
be the set of the first $k$ vertices numbered by $\eta$. We observe
that, in any interval graph, for any three vertices $u$, $v$ and $w$
such that $\eta(u)<\eta(v)<\eta(w)$, if $w$ is adjacent to $u$, then
$w$ is also adjacent to $v$ and therefore $N\left(S_k(\eta), I_{\min} \right)=\{u\in
N\left(\eta^{-1}(\nscomm{k})\right)| \eta(u)>k\}$. Using this fact, we have
\[
\sum_{k=1}^{n-1}|\vb\left(S_k(\eta),G\right)| \le 
\sum_{k=1}^{n-1}|\vb\left(S_k(\eta), I_{\min} \right)|=\sum_{u\in V}|\left\{v|
uv\in E\mbox{ and }\eta(v)>\eta(u)\right\}|=\numed{\imin}.
\]
By definition, $|\vb\left(S_k(\eta),G\right)|\ge b_v(k,G)$ and therefore,
$\numed{\imin}\ge\sum_{k=1}^{n-1}b_v(k,G)$. Now, the upper bound on
$\numed{\comp{\imin}}$ follows by applying equation (\ref{eqn:cvk}):
\begin{eqnarray*}
\numed{\comp{\imin}}&&={n\choose 2}-\numed{\imin}\\ 
&&\le{n\choose 2}-\sum_{k=1}^{n-1}b_v(k,G)=\sum_{k=1}^{n-1}n-k-b_v(k,G)
=\sum_{k=1}^{n-1}c_v(k,\comp{G}).  
\end{eqnarray*}
\bqed
\end{proof}
\subsection{Some Immediate Applications of Theorem \ref{thm:cvithm}}
\label{sec:immAppMethod1}
\begin{theorem}\label{thm:nk2}
Let $G(V,E)$ be a $(n-k-1)$-regular graph for some positive integer $k$.
Then, $\boxi(G)\ge\frac{n}{2k}$. 
\end{theorem}
\begin{proof} 
Since $\comp{G}$ is a $k$-regular graph, $\numed{\comp{G}}=nk/2$
and $c_v(1,\comp{G})=k$.  For $2\le i\le k$, $c_v(i,\comp{G})\le k$
(by Observation \ref{obs:cvObs1}) and for $i>k$, $c_v(i,\comp{G})=0$ (by
Observation \ref{obs:cvObs2}). Hence, $\sum_{i=1}^{n-1}c_v(i,\comp{G})\le
k^2$. Applying Theorem \ref{thm:cvithm}, the result follows.\bqed
\end{proof}
\paragraph{\textnormal{\textbf{Tightness of Theorem \ref{thm:nk2}:}}}
Consider the following co-bipartite graph $G(V,E)$: $V=A\uplus B$
and let $n$, $k$ and $l$ be integers such that, $|A|=|B|=kl=n/2$,
further, let $A=A_1\uplus A_2\uplus\cdots\uplus A_l$ and $B=B_1\uplus
B_2\uplus\cdots\uplus B_l$, where, $|A_i|=|B_i|=k$. For any two vertices,
$u$ and $v$, let $uv\notin E\Longleftrightarrow u\in A_i \mbox{ and }
v\in B_i$, $1\le i\le l$. Clearly, $G$ is an $(n-k-1)$-regular graph and
therefore by Theorem \ref{thm:nk2}, $\boxi(G)\ge n/2k=l$. Now, we present
a set of $l$ interval graphs $\mathcal{I}=\{I_i| 1\le i\le l\}$ whose edge
intersection gives $G$. The interval representation $f_i$ for each $I_i$
is as follows: $f_i(u)=[0,1]$ if $u\in A_i$, $f_i(u)=[2,3]$ if $u\in
B_i$, and $f_i(u)=[1,2]$ otherwise. Therefore, $\boxi(G)=l=n/2k$. In fact,
we note that each $I_i$ is a unit-interval graph by construction.
Therefore, this graph acts as a tight example even if we were to replace
$\boxi(G)$ with $\cub(G)$ in Theorem \ref{thm:nk2}.

\begin{theorem} \label{thm:nkExamples}
Let $k\ge 1$ be an integer.
\begin{enumerate}
\item Let $G$ be a ($n-k-1$)-regular co-planar graph. Then,
$\boxi(G)\ge\frac{n}{8}$.
\item Let $G$ be the complement of a $k$-regular $C_4$-free graph. Then,
$\boxi(G)\ge\frac{n}{4}$.
\item $\boxi(\comp{C_n})\ge\frac{n}{3}$.
\end{enumerate}
\end{theorem}
\begin{proof}
(1) Since $\comp{G}$ is a planar graph, it does not contain $K_{3,3}$
as a subgraph. Hence, $c_v(i,\comp{G})\le 2$ for $3\le i\le k$. From
Observation \ref{obs:cvObs2}, $c_v(i,\comp{G})=0$ for $i>k$. The result
follows by applying Theorem \ref{thm:cvithm} with $c_v(i,\comp{G})\le k$,
for $i=1,2$.\\
(2) Since $\comp{G}$ is $C_4$ or $K_{2,2}$-free, the result follows in the
same way as for (1).\\
(3) We have $\numed{C_n}=n$, $c_v(1,C_n)=2$, $c_v(2,C_n)=1$ and by
Observation \ref{obs:cvObs2}, $c_v(i,C_n)=0$ for $i>2$. Applying
Theorem \ref{thm:cvithm} the result follows. We recall that Cozzens
and Roberts \cite{computingBoxicityCozzensRoberts} showed that
$\boxi(\comp{C_n})=\ceil{n/3}$ by different methods. \bqed
\end{proof}
\begin{remark}
Theorems \ref{thm:nk2} and \ref{thm:nkExamples} are only indicative. We can
bound the boxicity of various other graph classes using this
method. If $G$ is not regular we can give a lower bound in terms of the
ratio $\delta(\comp{G})/\Delta(\comp{G})$, i.e. $\boxi(G)\ge
\frac{n\delta(\comp{G})}{2\Delta(\comp{G})^2}$, when 
$\Delta(\overline G) > 0$. Later, using the second method
we will show a better lower bound of $\frac{n}{2\Delta(\comp{G})}$, 
when $G$ does not contain universal vertices.
\end{remark}

\subsection{Lower Bounds for Random Graphs: $\mathcal{G}(n,p)$ and
$\mathcal{G}(n,m)$}\label{sec:rand1}
Following upper
bound on ${n \choose k}$, which is a consequence of Stirling's approximation for $n!$,  will be used frequently.
For  $1\le k\le n$,
\begin{equation} \label{eqn:nci} 
{n \choose
k}<\left(\frac{en}{k}\right)^k \le \exp\left(2k\log(n/k)\right).
\end{equation}


Also, we
use the following  variants of Chernoff's bound \cite {probCompMitzenmacherUpfal}. Let
$X$ be the sum of $n$ independent Poisson trials. Then 

$$\Pr[X \ge (1+\epsilon) E(X)] \le \exp(-\min\{\epsilon, \epsilon^2\}E(X)/3) \mbox{~~for any~~} \epsilon > 0;$$
$$\Pr[X \le (1-\epsilon) E(X)] \le \exp(-\epsilon^2 E(X)/2) \mbox{~for~} 0 < \epsilon < 1.$$


\ignore{
	\begin{theorem}\label{thm:randThm}
	Let $G\in \mathcal{G}(n,p)$ and $p$ be such that 
	$c_1/n\le
	p\le1-c_2\frac{\log n}{n^2}$, where $c_1$ and $c_2$
	are suitably chosen positive constants. Then, for almost all graphs
	$G\in \mathcal{G}(n,p)$, $\boxi(G)=\Omega(np(1-p))$.
	\end{theorem}
}

\begin{theorem}\label{thm:randThm}
Let $p\le1-\frac{ 40 \log n}{n^2}$.
Then
for $G\in \mathcal{G}(n,p)$, $\boxi(G)=\Omega(np(1-p))$ with probability at least $1 - 3/n^2$.
\end{theorem}
\begin{proof}

Consider the event $\numed{\comp{G}} < \frac{1}{2} {n \choose 2}(1-p)$.
We have $E(\numed{\comp{G}}) = {n \choose 2}(1-p)$.
By Chernoff's bound, we obtain that
$\Pr(\numed{\comp{G}} < E(\numed{\comp{G}})/2) \le \exp(- E(\numed{\comp{G}})/8) \le 1/n^2$.

In the rest of the proof, we assume that $p \ge c/n$ where $c$ is a suitably fixed constant. 
It is easy to see that when $p < c/n$, the  $\Omega(np(1-p))$ 
lower bound given by the theorem trivially holds true,
 by the definition of $\Omega$-notation.

We show that the event
$(\numed{\comp{G}} \ge \frac{1}{2} {n \choose 2}(1-p)$ and 
$\sum_{i=1}^{n-1}c_v(i,\comp{G}) \le 62n/p)$ holds true
with probability at least $1 - 3/n^2$.
The result then follows
directly from Theorem \ref{thm:cvithm}.

It suffices to prove that
\(
\Pr( \sum_{i=1}^{n-1}c_v(i,\comp{G}) > 62 n/p) \le 2/n^2.
\)
We divide the summation $\sum_{i=1}^{n-1}c_v(i,\comp{G})$
into following two parts.

\paragraph{\textbf{Range \textrm{1} $\left(1\le i\le\frac{5}{p}\log (np)\right)$:}} Let $k = \lfloor \frac{5}{p}\log (np) \rfloor$
and let $r_i=7n(1-p)^i + 6n/\log (np)$.
Note that if $c_v(i, \comp{G}) \le r_i$ for each $i \in \{1, \ldots, k\}$, then clearly $\sum_{i=1}^k c_v(i, \comp{G}) \le \sum_{i=1}^k r_i \le 37n/p$.
In the following we show that for $i \in \{1, \ldots, k\}$ we obtain $\Pr(c_v(i, \comp{G}) > r_i) \le 1/n^3$.
By union bound it follows that $\Pr(\sum_{i=1}^k c_v(i, \comp{G}) \le 37n/p) \ge 1 - 1/n^2$.

Consider any $A\subseteq V$ with $|A|=i$. Consider the random
variable $X=|\sn(A,\comp{G})|$ and the event $X > r_i$.
Note that $E(X) = (n-i)(1-p)^i$. By Chernoff's bound, 
we obtain
\(
Pr(X > r_i)  \le Pr \left (X > (1 + \frac {r_i -E(X)} {E(X)} ) E(X) \right ) \le \exp(-2n(1-p)^i -2n/\log(np)) \le \exp(-2n/\log(np)).
\)
Recalling the definition of $c_v(i, \comp{G})$, it follows that $\Pr(c_v(i, \comp{G}) > r_i)  \le {n \choose i} \exp(-2n/\log(np)) \le {n \choose k} \exp(-2n/\log(np))$ as $k \le \frac{5}{p} \log(np) \le n/2$ when $p \ge c/n$ for a suitably chosen constant $c$.
Using (\ref{eqn:nci}) and the fact that $k = \lfloor \frac{5}{p} \log(np) \rfloor$, it is straightforward to verify that ${n \choose k} \exp(-2n/\log(np)) \le \exp(-(2n/\log(np) - 2k \log(n/k))) \le \exp(-n/\log(np))$ when $p \ge c/n$ for suitably fixed constant $c$. Clearly $\exp(-n/\log(np)) = o(1/n^3)$.

\paragraph{\textbf{Range \textrm{2} $\left(\frac{5}{p}\log (np) \le i \le n\right)$:}} 
Let $k =\left  \lceil \frac{5}{p}\log (np) \right \rceil$.
Note that if $c_v(i, \comp{G}) \le \frac{10}{p} \log(4n/i)$ for each $i \in \{k, \ldots, n\}$, 
then $\sum_{i=k}^n c_v(i, \comp{G}) \le \sum_{i=k}^n \frac{10}{p}\log(4n/i) \le \frac{10}{p} \sum_{i=1}^n \log(4n/i) \le 25n/p$.
The last inequality follows easily from Stirling's inequality that $\log(n!) \ge n \log(n) - n$.
In the following we show that for $i \in \{k, \ldots, n\}$ we obtain $\Pr(c_v(i, \comp{G}) \ge \frac{10}{p}\log(4n/i)) \le 1/n^3$.
By union bound it follows that $\Pr(\sum_{i=k}^n c_v(i, \comp{G}) \le 25n/p) \ge 1 - 1/n^2$.

For any two disjoint sets $A, B \subseteq V$ with $|A| = x$ and $|B| = y$,
we have $\Pr(B \subseteq N_S(A, \comp{G})) = (1-p)^{xy}$. Since $c_v(x, \comp{G}) \ge y$ implies
existence of such $A$ and $B$, by union bound it follows that 
\begin{equation}\label{eqn:count}
\Pr(c_v(x, \comp{G}) \ge y) \le {n \choose x}{n \choose y} (1-p)^{xy}
\end{equation}
Using (\ref{eqn:nci}) and the fact that $(1-p) \le e^{-p}$, we have 
$\Pr(c_v(x, \comp{G}) \ge y) \le \exp(-(pxy - 2x \log(n/x) - 2y \log(n/y)))$.
Choosing $x = i$ and $y = \frac{10}{p} \log(4n/i)$  it follows that 
$\Pr(c_v(i, \comp{G}) \ge \frac{10}{p} \log(4n/i)) \le \exp(-(10i \log(4n/i) - 2i \log(n/i) - \frac{20}{p} \log(4n/i) \log(np))) \le \exp(-(8i \log(4n/i) - \frac{20}{p} \log(4n/i)\log(np))) \le \exp(-4i \log(4n/i))$ where the last inequality follows from the fact that $k \ge \frac{5}{p} \log(np)$.  
It is straightforward to verify that $\exp(-4i \log(4n/i)) = (i/4n)^{4i} =  o(1/n^3)$.

From above two cases we conclude that $\Pr(\sum_{i=1}^n c_v(i, \comp{G}) \le 62n/p) \ge 1 - 2/n^2$.\bqed
\end{proof}

\begin{corollary}\label{cor:npRandThm}
Let $G\in \mathcal{G}(n,p)$. Let $p$ be such that $p\le c<1$, where
$c$ is any positive constant. Then  $\boxi(G)=\Omega(np)$ with probability at least $1 - 3/n^2$.
\end{corollary}
\begin{corollary}\label{cor:randThmEdges}
Consider graphs with $n$ vertices and $m$ edges such that $m\le c n^2/3$, where $c < 1$ is any positive constant.
Then for almost all graphs $G$ on $n$ vertices and $m$ edges, $\boxi(G) = \Omega(m/n)$.
\end{corollary}
\begin{proof}
Let $N = {n \choose 2}$ and let $p =\frac{m}{N}$.
It can be shown that, given any property $\mathcal{P}$,
$Pr(G \in \mathcal{G}(n,m) \mbox{ does not satisfy }\mathcal{P})\le
3\sqrt{m}Pr(G \in \mathcal{G}(n, p) \mbox{ does not satisfy }\mathcal{P})$
(see \cite[Theorem 2 Chapter 2]{bollobasRandomGraphs}).
Let $\mathcal{P}$ be the property that $\boxi(G) = \Omega(np)$.
From Corollary \ref{cor:npRandThm}, we have $\Pr(\boxi(G) = \Omega(np)) \ge 1 - 3/n^2$,
when $p \le c$. Recalling that $p = m/N$, it follows that
$Pr(\mbox{For~} G \in \mathcal{G}(n,m), \boxi(G) = \Omega(nm/N))\ge 1 - 9 \sqrt{m}/n^2$
when $m/N \le c$. In other words, for $G \in \mathcal{G}(n,m)$, $\Pr(\boxi(G) = \Omega(m/n)) \ge 1 - 9/n$
when $m \le cn(n-1)/2$. Hence proved. \bqed
\end{proof}

\ignore{
	We recall that, in the proof of Theorem \ref{thm:randThm},
	$\boxi(G)=\Omega(np(1-p))$ if the event $Q$ does not occur. Let
	$p=\frac{m}{{n\choose 2}}$. 
	It can be shown that, given any property $\mathcal{P}$,
	$Pr(G \mbox{ satisfies }\mathcal{P}, G\in\mathcal{G}(n,m))\le
	3\sqrt{m}Pr(G \mbox{ satisfies }\mathcal{P},
	G\in\mathcal{G}(n,p))$ (see \cite[Theorem 2 Chapter 2]{bollobasRandomGraphs}).
	Hence, $Pr(G\mbox{ satisfies } Q,G\in\mathcal{G}(n,m))\le
	3\sqrt{m}\mathcal{O}(1/n^2)=\mathcal{O}(1/n)$. Therefore, for almost all
	$G\in\mathcal{G}(n,m)$, $\boxi(G)=\Omega(m/n)$. Hence proved.\bqed
}

\ignore{
	For very dense graphs, i.e. with $p\rightarrow1$ as $n\rightarrow\infty$,
	it is not possible to infer from Theorem \ref{thm:randThm} that
	$\boxi(G)=\Omega(np)$ because of the additional factor $(1-p)$.
	For this case we use a different approach to bound the boxicity. A
	theorem follows in this regard.
}

In the following, we derive another lower bound for boxicity of $G \in \mathcal{G}(n, p)$.
This bound is interesting in the sense that for very dense random graphs, 
that is for the cases when $p \rightarrow 1$ as $n \rightarrow \infty$,
it  is stronger 
 than $\Omega(np(1-p))$ bound given by Theorem \ref{thm:randThm}.

\begin{lemma}\label{lem:randlow2}
Let $G \in \mathcal{G}(n, p)$. Then  $\sum_{i=1}^{n-1} c_v(i, \comp{G})= O\left(\frac{n(1-p)\log n + \log^2(n)}{\log(1/(1-p))} \right)$ 
with probability at least $1-2/n^2$.
\end{lemma}
\begin{proof}

Let $\bar{d}(u)$ denote the degree of vertex $u$ in $\comp{G}$. We have $E(\bar{d}(u)) = (n-1)(1-p)$.
By Chernoff's bound, we obtain $\Pr(\bar{d}(u) \ge 10 (\log n + n(1-p)) ) \le  \exp(- 3 (\log n + n (1-p))) \le 1/n^3$.
By union bound, it follows that 
\begin{equation}\label{eqn:cv1}
\Pr(c_v(1, \comp{G}) \ge 10  (\log n + n (1-p))) \le 1/n^2 
\end{equation}

Let $r = \left  \lceil\frac{5\log(n)}{\log(1/(1-p))} \right \rceil$. 
For any fixed $i > r$, we now bound 
$\Pr(c_v(i, \comp{G}) \ge r)$.
Using (\ref{eqn:count})  and the fact that ${n \choose i} \le n^i$, we  obtain
\begin{equation}\label{eqn:cvj}
Pr\left(c_v(i,\comp{G})\ge r\right) \le \exp\left(-(ir\log(1/(1-p)) - i\log n - r \log n )\right) \le \exp(-3\log n) = 1/n^3
\end{equation}
where the last inequality follows from the fact that $i > r$ and that $r \ge \frac{5 \log n}{\log(1/(1-p))}$.

From Observations \ref{obs:cvObs1} and \ref{obs:cvObs2},
it is straightforward to see that
$$
\sum_{i=1}^{n-1}c_v(i,\comp{G}) =\sum_{i=1}^{c_v(1,\comp{G})}c_v(i,\comp{G}) \le \sum_{i=1}^{r}c_v(1,\comp{G}) + \sum_{i=r+1}^{c_v(1,\comp{G})}c_v(i,\comp{G}) =  r \cdot c_v(1,\comp{G}) + \sum_{i=r+1}^{c_v(1,\comp{G})}c_v(i,\comp{G}).
$$

From (\ref{eqn:cv1}) and (\ref{eqn:cvj}),
above equation implies the result that
$\sum_{i=1}^{n-1}c_v(i,\comp{G}) = O(r \cdot c_v(1, \comp{G})) = O\left(\frac{n(1-p)\log n + \log^2(n)}{\log(1/1-p)}\right)$
with probability at least $1  - 2/n^2$.  \bqed
\end{proof}

\begin{theorem}\label{thm:dense}
Let $p \le 1 -  \frac{40 \log n}{n^2}$.
Let $\gamma = \frac{n(1-p)\log n + \log^2(n)}{\log(1/(1-p))}$.
Then for $G \in \mathcal{G}(n, p)$, 
$\boxi(G) =\Omega(n^2 (1-p)/\gamma)$ with probability at least $1 - 3/n^2$.
\end{theorem}
\ignore{
	\begin{theorem}\label{thm:dense}
	Let $G\in \mathcal{G}(n,p)$ with $p=1-g/n$, where $g(n)$ is any
	non-negative function such that $g(n)=o(n)$ and $g(n)=\omega(1/n)$. Then,
	\begin{enumerate}
	\item if $g\ge c\log n$ where $c$ is some constant, then for almost all
	graphs $\boxi(G)=\Omega\left(n\right)$,
	\item if $g<c\log n$, then for almost all graphs,
	$\boxi(G)=\Omega\left(\frac{ng}{\log n}\right)$,
	\end{enumerate}
	\end{theorem}
}
\begin{proof}
We have $E(\numed{\comp{G}}) = {n \choose 2}(1-p)$.
By Chernoff's bound, we obtain that
$\Pr(\numed{\comp{G}} < E(\numed{\comp{G}})/2) \le \exp(- E(\numed{\comp{G}})/8) \le 1/n^2$.
The result now follows
directly from Theorem \ref{thm:cvithm} and Lemma \ref{lem:randlow2}.  \bqed
\end{proof}

\begin{corollary}\label{cor:npDenseRandThm1}
Let $G \in \mathcal{G}(n, p)$ and let $1 - \frac{\log n}{n} \le p \le 1 -
 \frac{40 \log n}{n^2}$.
Then  $\boxi(G) =\Omega(n^2 (1-p)/\log n))$ with probability at least $1 - 3/n^2$.
\end{corollary}

\begin{corollary}\label{cor:npDenseRandThm2}
Let $G \in \mathcal{G}(n, p)$ and let  $p < 1 - \log n/n$.
Then $\boxi(G) =\Omega(\frac{n}{\log n}\log(\frac{1}{1-p}))$ with probability at least $1 - 3/n^2$.
\end{corollary}

\ignore{
	We note that the above bound is weaker than Theorem \ref{thm:randThm} bound for small values of $p$. Already when $p$ is a constant (say $1/2$),
	the later bound (Corollary \ref{cor:npDenseRandThm2}) only yields $\Omega(n/\log n)$ bound. For $p < 1/2$,
	we have $1-p > p$. Also, $1/(1-p) = 1 + (p/(1-p)) \le \exp(p/(1-p))$, since $p/(1-p) < 1$. Hence
	$\log(1/(1-p)) \le p/(1-p) \le 2p$ since $1-p > 1/2$. Thus, for $p <1/2$, above corollary yields $\Omega(\frac{n}{\log n} \log(1/(1-p)))$ which 
	is no better than $\Omega(np/\log n)$ which is weaker than $\Omega(np)$ lower bound of Theorem \ref{thm:randThm}.
}

We note that for $1 - \frac{\log\log n}{\log n} \le p \le 1 -  \frac{40\log n}{n^2}$, Corollaries \ref{cor:npDenseRandThm1} and \ref{cor:npDenseRandThm2} give improved lower bounds than $\Omega(np(1-p))$ lower bound given by Theorem \ref{thm:randThm}. In particular, 
it follows that for almost all graphs $G \in \mathcal{G}(n, p)$, 
$\boxi(G) = \Omega(n)$ when $p =1-  \frac{\log n}{n}$ and $\boxi(G) = \Omega(n/\log n)$ when $1 - 1/\log n \le p \le 1 - 1/n$.
\ignore{
	Following corollary is a direct consequence of above two corollaries.
	\begin{corollary}\label{cor:npDenseRandThm3}
	Let $G \in \mathcal{G}(n, p)$ and let  $1 - 1/\log n \le p \le 1 - 1/n$.
	Then $\boxi(G) =\Omega(n/\log n)$ with probability at least $1 - 3/n^2$.
	\end{corollary}
}

\subsection{Spectral Lower Bounds}\label{sec:spec}
Consider a graph $G(V,E)$ and let  $X\subseteq V$.  Recalling the
definition of vertex-boundary and vertex neighbourhood from Section
\ref{sec:neighbours}, we have $\vb(X)\subseteq N'(X)$ and
\begin{equation}
\label{eqn:phiN} |\vb(X)| \ge |N'(X)|-|X|.
\end{equation}

Suppose $G(V,E)$ is a $k$-regular balanced $(A,B)$-bipartite graph with
$2n$ vertices. For any $X\subseteq A$, Tanner \cite{concentratorsTanner}
gives a lower bound for $|N'(X)|/|X|$ using spectral methods. Suppose
the vertices in $A$ and $B$ are ordered separately, we define
the \emph{bipartite incidence matrix} $M$ as follows: $M=[m_{ij}]$,
$m_{ij}=1$ if the $i$th vertex in $A$ is adjacent to $j$th vertex in
$B$. The largest eigenvalue of $MM^T$ is $k$. Let $\lambda'$ be the
second largest eigenvalue of $MM^T$. Then,
\begin{equation} \label{eqn:tannerBipartite}
|N'(X)| \ge \frac{k^2|X|}{\lambda'+(k^2-\lambda')\frac{2|X|}{n}}.
\end{equation}
This result can be extended to any $k$-regular graph $G(V,E)$ (see
\cite{phdThesisNabil}) as follows. Let $\lambda$ be the second largest
eigenvalue in absolute value of the adjacency matrix of $G$. For any
subset $X\subseteq V$,
\begin{equation} \label{eqn:tanner}
|N'(X)| \ge \frac{k^2|X|}{\lambda^2+(k^2-\lambda^2)\frac{|X|}{n}}.
\end{equation}
\begin{theorem}\label{thm:spec}
Let $G(V,E)$ be a connected $k$-regular graph on $n$ vertices. Let
$\lambda$ be the second largest eigenvalue in absolute value of the
adjacency matrix of $G$. Then,
\[
\boxi(G)\ge\left(\frac{k^2/\lambda^2}{\log\left(\frac{k^2}{\lambda^2}+1\right)}\right)\left(\frac{n-k-1}{2n}\right).
\]
\end{theorem}
\begin{proof}
We first note that $\numed{\comp{G}}=\frac{n(n-k-1)}{2}$. By the
definition of $b_v(i)$ in equation (\ref{eqn:bvi}) (See Section
\ref{sec:firstMethod}) and by applying equations (\ref{eqn:phiN}) and
(\ref{eqn:tanner}):
\begin{equation}\label{eqn:bviub}
b_v(i,G) \ge \frac{k^2i}{\lambda^2+(k^2-\lambda^2)\frac{i}{n}} - i.
\end{equation}
Recalling that $c_v(i,\comp{G})=n-i-b_v(i,G)$ and applying Theorem
\ref{thm:cvithm} we get,
\begin{eqnarray*}
\boxi(G) &\ge&  
\frac{n(n-k-1)/2}{\ds\sum_{i=1}^{n-1}n-i-\left(\frac{k^2i}{\lambda^2+(k^2-\lambda^2)\frac{i}{n}}
- i\right)} \\
&=&
\frac{(n-k-1)/2}{\ds\sum_{i=1}^{n-1}1-\frac{k^2i}{\lambda^2n+(k^2-\lambda^2)i}
} \\
&\ge&
\frac{(n-k-1)/2}{\ds\sum_{i=1}^{n-1}1-\frac{k^2i}{\lambda^2n+k^2i}} =
\frac{(n-k-1)/2}{\ds\sum_{i=1}^{n-1}\frac{\lambda^2n}{\lambda^2n+k^2i}} \\
&=&
\frac{(n-k-1)/2}{\ds\frac{\lambda^2n}{k^2}\ds\sum_{i=1}^{n-1}\frac{1}{\frac{\lambda^2n}{k^2}+i}} \\
&\ge&
\frac{(n-k-1)/2}{\ds\frac{\lambda^2n}{k^2}\left(\log\left(\frac{\lambda^2n}{k^2}+n\right)-\log\left(\frac{\lambda^2n}{k^2}\right)\right)} \\
&=&
\frac{k^2/\lambda^2}{\log\left(\frac{k^2}{\lambda^2}+1\right)}\left(\frac{n-k-1}{2n}\right).
\end{eqnarray*}
Hence proved. \bqed
\end{proof}
\subsubsection{Tightness of Theorem \ref{thm:spec}}
Let $l$ and $p$ be integers $\ge 2$. Consider a complete $p$-partite graph
$K_{l,l,\ldots,l}$ where $n=lp$. $K_{l,l,\ldots,l}$ is a \emph{strongly
regular graph}, i.e. it is a $k$-regular graph where each pair of
adjacent vertices has the same number $a\ge0$ of common neighbours,
and each pair of non-adjacent vertices has the same number $c\ge1$ of
common neighbours. For a strongly regular graph, $\lambda$ can be obtained
by solving the quadratic equation $x^2+(c-a)x+(c-k)$ (See \cite[Chapter
3]{algGraphTheoryBiggs}). For $K_{l,l,\ldots,l}$, noting that $k=c=(p-1)l$
and $a=(p-2)l$, we get $\lambda=l$. Using Theorem \ref{thm:spec}
we get $\boxi(G)=\Omega\left(\frac{p}{\log p}\right)$. In
\cite{recentProgressesInCombRoberts}, it is
shown that $\boxi(K_{l,l,\ldots,l})=p$. Hence,
our result is tight up to a $\mathcal{O}(\log(p)) =
\mathcal{O}\left(\log\left(\frac{k^2}{\lambda^2}+1\right)\right)$
factor. It would be interesting to see if we can improve this bound to
$\boxi(G)=\Omega\left(\frac{k^2}{\lambda^2} \left(\frac{n-k-1}{2n}\right)\right)$. For
example, when $l=2$ (This graph was considered by Roberts in
\cite{recentProgressesInCombRoberts} and can be obtained by
removing a perfect matching from a complete graph), such an improved
bound will evaluate to $\frac{n-2}{8}$. The boxicity of this graph
was shown to be $n/2$ in \cite{recentProgressesInCombRoberts}.
\subsubsection{Applying the Spectral Bound to Random Regular Graphs}
In 1986, Alon and Boppana \cite{eigenvaluesExpandersAlon} proved that
for a fixed $k$, $\lambda\ge2\sqrt{k-1}$ for almost all $k$-regular
graphs. In the same paper it was conjectured that for any $k\ge3$ and
$\epsilon>0$, $\lambda\le2\sqrt{k-1}+\epsilon$. Recently, Friedman
\cite{proofAlonSecondEigenvalueConjectureFriedman04} proved this
conjecture albeit for a different model of regular graphs (involving
multiple edges and loops). However, using contiguity theorems
\cite{modelsRandomGraphsWormald} it is easy to infer that the same
result applies to the $\mathcal{G}_R(n,k)$ model too. Hence, we can
assume that $\lambda\approx2\sqrt{k-1}$ for a random $k$-regular graph.
We have the following corollary to Theorem \ref{thm:spec}.
\begin{corollary}\label{cor:specAlmost}
Let $\mathcal{G}_R(n,k)$ be the probability space of random $k$-regular
graphs, where $k$ is fixed. For  almost all $G\in\mathcal{G}_R(n,k)$,
$\boxi(G)=\Omega\left(k/\log k\right)$.
\end{corollary}

\subsection{Lower Bounds for Bipartite Graphs}
The presence of large independent sets or large sparse induced subgraphs
 in $G$ puts a limitation on the
lower bound that can be achieved using Theorem \ref{thm:cvithm}. For
example, consider a bipartite graph on $2n$ vertices. It contains an
independent set of size at least $n$, the best possible bound we can
get by applying Theorem \ref{thm:cvithm} is ${2n\choose 2}/{n\choose 2}
\le 5$, for $n \ge 3$. But, by using a simple trick we can overcome
this limitation to a great extent. Our basic tool is the following lemma
from \cite{constFactorApproxCirArcAdigaBabuChandran}.

\begin{lemma}\label{lem:ind2Cliq}
Let $G(V,E)$ be a graph with a partition $(V_1,V_2)$ of its vertex set $V$.
Let $G_1(V,E_1)$ be its supergraph such that $E_1 = E \cup \{ (a',b'):
a',b' \in V_2 \}$. Then $\boxi (G_1) \le 2 \boxi(G)$.
\end{lemma} 

If $G$ has a large independent set (or a large sparse induced subgraph),
 the above lemma allows us to get a
new graph $G_1$, where the independent set is replaced by a clique on
the same vertex set and everything else is as in $G$. Now, we can try to
get a lower bound for $\boxi(G_1)$. Clearly, the  lower bound for $\boxi(G_1)$
multiplied by $1/2$ will give a lower bound for $\boxi(G)$. Note that
if there are more than one disjoint independent sets in $G$, 
we can apply this trick multiple times. 

In this section
we will derive lower bounds (that correspond to Theorems \ref{thm:nk2},
\ref{thm:randThm}, \ref{thm:dense} and \ref{thm:spec}) for balanced bipartite
graphs, using this method in conjunction with Theorem \ref{thm:cvithm}.
Let $G=(V,E)$ be a bipartite graph with bipartition $(V_1,V_2)$.
Let $G'=(V,E')$ be the graph with edge set $E' = E \cup \{ (a',b'):
a',b' \in V_2 \} \cup \{ (a',b'): a',b' \in V_1 \}$. By applying Lemma
\ref{lem:ind2Cliq}, we infer that $\boxi(G')\le4\boxi(G)$. However,
we can improve the bound to $2\boxi(G)$ by slightly modifying 
Lemma~\ref{lem:ind2Cliq}.

\begin{lemma}\label{lem:bip2Cobip}
Let $G(V,E)$ be a graph with a partition $(V_1,V_2)$ of its vertex set $V$.
Let $G'(V,E')$ be its supergraph such that $E' = E \cup \{ (a',b'):
a',b' \in V_1 \} \cup  \{ (a',b'): a',b' \in V_2 \}$.  Then, $\boxi(G')
\le 2 \boxi(G)$.
\end{lemma} 
\begin{proof}
Let $\boxi(G)=b$ and let $\igr=\{I_1,I_2,\ldots,I_b\}$ be an interval graph
representation of $G$. For each $I\in\igr$, we construct two interval
graphs $I^1$ and $I^2$ such that $\igr'=\{I_i^j|\ j=1,2,\ 1\le i\le b\}$
is an interval graph representation of $G'$. Let $f$ be an interval
representation of $I$. Without loss of generality we will assume that for all
$u\in V$, $0\le l_f(u),r_f(u)\le n$. Let $f^1$ correspond to an interval
representation of $I^1$ which is realized by modifying $f$ as follows: If
$u\in V_1$, then, $l_{f^1}(u)=0$ and $r_{f^1}(u)=r_f(u)$ and if $u\in V_2$,
then, $r_{f^1}(u)=n$ and $l_{f^1}(u)=l_f(u)$. Similarly, the interval
representation for $I^2$, $f^2$ is realized as follows: If
$u\in V_2$, then, $l_{f^2}(u)=0$ and $r_{f^2}(u)=r_f(u)$ and if $u\in V_1$,
then, $r_{f^2}(u)=n$ and $l_{f^2}(u)=l_f(u)$. It is straightforward to
verify that $\igr'$ is indeed an interval graph representation of $G'$.
\bqed
\end{proof}

Now in order to derive lower bounds for the boxicity of a bipartite
graph $G$, we can use $G'$ instead. 

\begin{theorem}\label{thm:bipnk2}
Let $G$ be a $(n-k)$-regular balanced bipartite graph on $2n$
vertices. Then, $\boxi(G)\ge\frac{n}{2k}$.
\end{theorem}
\begin{proof}
Construct $G'$ from $G$ as in Lemma \ref{lem:bip2Cobip}. Now $G'$ is
a $(2n-k-1)$-regular graph on $2n$ vertices. By Theorem \ref{thm:nk2},
$\boxi(G')\ge\frac{n}{k}$. Applying Lemma \ref{lem:bip2Cobip}, we get $\boxi(G)
\ge\frac{n}{2k}$. \bqed
\end{proof}

Considering the rather indirect trick used to obtain the lower bound
of Theorem \ref {thm:bipnk2}, one may think that it may not be tight. But
surprisingly, it turns out to be a tight lower bound!

\paragraph{\bf Tightness of Theorem \ref{thm:bipnk2}:} We now give a
tight example for the above theorem. Let $G$ be a balanced $(V_1,V_2)$-bipartite
graph on $2n$ vertices with the following structure: Let $n$, $k$ and $l$
be integers such that $l$ is even and $kl=n$. Further, let $V_1=A_1\uplus
A_2\uplus\cdots\uplus A_l$ and $V_2=B_1\uplus B_2\uplus\cdots\uplus
B_l$, where, $|A_i|=|B_i|=k$. For any two vertices, $u\in A_i$ and
$v\in B_j$, let $uv\in E\Longleftrightarrow i\ne j$. Clearly, $G$ is an
$(n-k)$-regular graph. From Theorem  \ref{thm:bipnk2}, $\boxi(G)\ge
n/2k=l/2$. Now, we present an interval graph representation of $G$,
$\mathcal{I}=\{I_i| 1\le i\le l/2\}$ using $l/2$ interval graphs.

The first interval graph $I_1$ is defined as follows: Each vertex $v\in
V_1$ is assigned a unique point interval $[i,i]$, $i\in\{1,2,\ldots,n\}$
such that if $v\in A_1$, then, $i\in\{1,2,\ldots,k\}$; if $v\in
V_1\setminus(A_1\cup A_2)$ then, $i\in\{k+1,k+2,\ldots,n-k\}$ and if $v\in
A_2$, then, $i\in\{n-k+1,k+2,\ldots,n\}$. Every $v\in B_2$ is assigned
$\left[1,n-k\right]$; $v\in B_1$ is assigned $\left[k+1,n\right]$ and
$v\in V_2\setminus(B_1\cup B_2)$ is assigned $\left[1,n\right]$.
The second interval graph $I_2$ is obtained
as follows: Each vertex $v \in V_2$ is assigned a unique point interval
$[i,i]$,   $i\in\{1,2,\ldots,n\}$
such that if $v\in B_3$, then, $i\in\{1,2,\ldots,k\}$; if $v\in
V_2\setminus(B_3\cup B_4)$ then, $i\in\{k+1,k+2,\ldots,n-k\}$ and if $v\in
B_4$, then, $i\in\{n-k+1,k+2,\ldots,n\}$. Every $v\in A_4$ is assigned
$\left[1,n-k\right]$; $v\in A_3$ is assigned $\left[k+1,n\right]$ and
$v\in V_1\setminus(A_3\cup A_4)$ is assigned $\left[1,n\right]$.
The other
interval graphs $I_j$, $j=3,\ldots,l/2$, are defined as follows: $v\in
A_{2j-1}$ is assigned $\left[1,2\right]$; $v\in B_{2j}$ is assigned
$\left[2,3\right]$; $v\in B_{2j-1}$ is assigned $\left[4,5\right]$; $v\in
A_{2j}$ is assigned $\left[5,6\right]$ and the rest of the vertices are
assigned the interval $[1,6]$.

It is easy to verify that each $I_j$ is an interval supergraph of $G$. In
$I_1$, $V_1$ induces an independent set and in $I_2$, $V_2$ induces an
independent set. For every pair $\{u,v\}$ such that $u\in A_k$, $v\in B_k$
$uv\notin E\left(I_{\ceil{k/2}}\right)$. Hence, $\igr$ is an interval graph
representation of $G$.

\subsubsection{Bounds for Random Bipartite Graphs} \label{sec:randbipartite}
\begin{theorem}\label{thm:bipartiteRandThm}
Let $p\le1-\frac{ 20 \log n}{n^2}$.  Then, for $G\in \mathcal{G}_B(2n,p)$,
$\boxi(G)=\Omega(np(1-p))$ with probability at least $1 - 3/n^2$.
\end{theorem}
\begin{proof}
Let $G\in \mathcal{G}_B(2n,p)$ and $(V_1,V_2)$ be the bipartition with
$|V_1|=|V_2|=n$. Let $G'$ be the corresponding cobipartite graph as
defined in Lemma \ref{lem:bip2Cobip}. We will show that for almost all
graphs $G\in \mathcal{G}_B(2n,p)$, $\boxi(G')=\Omega(np(1-p))$. From
Lemma \ref{lem:bip2Cobip}, it immediately follows that
$\boxi(G)=\Omega(np(1-p))$. The proof of the lower bound for $\boxi(G')$
is similar to that of Theorem \ref{thm:randThm}.

First, we observe that $\comp{G'}$ is a bipartite graph and in particular
$\comp{G'}\in \mathcal{G}_B(2n,1-p)$. Consider the event $\numed{\comp{G'}}
< \frac{1}{2}n^2(1-p)$. We have $E(\numed{\comp{G'}}) = n^2(1-p)$.
By Chernoff's bound, we have
$\Pr(\numed{\comp{G'}} < E(\numed{\comp{G'}})/2) \le \exp(-
E(\numed{\comp{G'}})/8) \le 1/n^2$.

We have the following observation to make:
\begin{observation}\label{obs:bipcv}
For a bipartite graph $H$ with bipartition $(V_1,V_2)$, if a set $Y$
has vertices from both $V_1$ and $V_2$, then, it will not have any strong
neighbor. Therefore, $c_v(\cdot)$ can be redefined as follows for a
bipartite graph:
\begin{equation}\label{eqn:bipcv}
c_v(k,H)=\max_{\stackrel{Y\subseteq V_1\mbox{\tiny or }V_2}{|Y|=k}}|\sn(Y,H)|.
\end{equation}
Moreover, when $Y \subseteq V_1$, all the strong neighbours of $Y$
will be in $V_2$ and vice versa.
\end{observation}
From the above observation, it also follows that $c_v(i,\comp{G'})=0$
for $i>n$. Now, the proof is almost identical to that of Theorem
\ref{thm:randThm} except for the following differences: (1) $n$ should
be interpreted as $|V_1|$ or $|V_2|$ instead of as $|V|$ ; (2) the sets $A$
and $B$ should be interpreted as subsets of $V_1$ and $V_2$ respectively
instead of as subsets of $V$ and (3) all occurrences of $(n-i)$ should
be replaced by $n$. \bqed
\end{proof} 

\begin{corollary}\label{cor:npBipartiteRandThm}
Let $G\in \mathcal{G}_B(n,p)$. Let $p$ be such that $p\le c<1$, where $c$
is any positive constant. Then  $\boxi(G)=\Omega(np)$ with probability at
least $1 - 3/n^2$.
\end{corollary}
\begin{corollary}\label{cor:randBipartiteThmEdges}
Consider balanced bipartite graphs with $2n$ vertices and $m$ edges such
that $m\le c n^2/3$, where $c < 1$ is any positive constant.  Then for
almost all balanced bipartite graphs $G$ on $2n$ vertices and $m$ edges,
$\boxi(G) = \Omega(m/n)$.
\end{corollary}
The proof is similar to that of Corollary \ref{cor:randThmEdges}. The
following results are analogous to Theorem \ref{thm:dense} and Corollaries
\ref{cor:npDenseRandThm1} and \ref{cor:npDenseRandThm2}.
\begin{theorem}\label{thm:bipariteDense}
Let $p \le 1 -  \frac{20 \log n}{n^2}$.  Let $\gamma = \frac{n(1-p)\log
n + \log^2(n)}{\log(1/(1-p))}$.  Then for $G \in \mathcal{G}_B(2n, p)$,
$\boxi(G) =\Omega(n^2 (1-p)/\gamma)$ with probability at least $1 -
3/n^2$.
\end{theorem}
\begin{corollary}\label{cor:npBipartiteDenseRandThm1}
Let $G \in \mathcal{G}_B(2n, p)$ and let $1 - \frac{\log n}{n} \le p \le
1 - \frac{20 \log n}{n^2}$. Then  $\boxi(G) =\Omega(n^2 (1-p)/\log n))$
with probability at least $1 - 3/n^2$.
\end{corollary}
\begin{corollary}\label{cor:npBipartiteDenseRandThm2}
Let $G \in \mathcal{G}_B(2n, p)$ and let  $p < 1 - \log n/n$.  Then $\boxi(G)
=\Omega(\frac{n}{\log n}\log(\frac{1}{1-p}))$ with probability at least
$1 - 3/n^2$.
\end{corollary}


\subsubsection{Spectral Bound}
\begin{theorem}
Let $G$ be a random $k$-regular balanced ($V_1,V_2$)-bipartite graph
on $2n$ vertices. Let $\lambda'$ be the second largest eigenvalue
of $MM^T$, where $M$ is the bipartite incidence matrix of $G$ (See
Section \ref{sec:spec} for definition). If $\lambda'\ne 0$, then,
\[
\boxi(G)\ge\left(\frac{k^2/\lambda'^2}{\log\left(\frac{k^2}{\lambda'^2}+1\right)}\right)\left(\frac{n-k}{2n}\right).
\]
\end {theorem} 
\begin {proof}
We construct $G'$ from $G$ as in Lemma \ref{lem:bip2Cobip}. The proof of
Theorem \ref{thm:spec} can now be applied to $G'$ with the
following minor modifications: (1) $||\comp{G'} || = n(n-k)$ instead of
$\frac{n(n-k-1)}{2}$, and (2) $c_v(i,\overline {G'}) = n - \frac{k^2
i}{\lambda'^2 + (k^2 - \lambda'^2)\frac{i}{n}}$. The lower bound for
$\boxi(G)$ is obtained by using Lemma \ref{lem:bip2Cobip}.
\bqed
\end {proof} 


\section{The Second Method}

\subsection{Lower Bound Theorem}

\begin{definition} \label{def_cross_exp}
Let $A,B \subseteq V$.  Let $t$ be a  positive integer such that
$t \le |A|$.  Then  $n_t (A,B)  = \min_{S \subseteq A, |S| = t} \{ |N[S,G] \cap
B| \}$. 
\end{definition}

\begin{definition} \label{def_comp_exp}
Let $A,B \subseteq V$. Let $t$ be a  positive
integer such that $t \le |A|$.  Then  let $m_t (A,B) 
= \min_{S
\subseteq A, |S| = t} \{ |N'(S,\comp{G}) \cap B| \}$.
\end{definition}

The following Lemma on interval graphs will be used in the 
proof of Theorem \ref {thm:thm_main}.
\begin{lemma} \label{lem_basic}
Let $I$ be an interval graph. Two vertices $u$ and $v$ are adjacent in $I$
if and only if $l(u) < r(v)$ and $l(v) < r(u)$.  Alternately, $u$ and $v$
are non adjacent in $I$ if and only if $r(u) < l(v)$ or $r(v) < l(u)$.
\end{lemma}
\newcommand{\dgone}[1]{\ensuremath{d_{(S_1,S_2)}(#1)}}
\newcommand{\dgtwo}[1]{\ensuremath{d_{(S_2,S_1)}(#1)}}
\newcommand{\dgmin}{\ensuremath{\delta_{(S_1,S_2)}}}
\newcommand{\dgmax}{\ensuremath{\Delta_{(S_2,S_1)}}}
\newcommand{\opt}{\ensuremath{b^*}}

\begin{theorem} \label{thm:thm_main}
Let $S_1, S_2 \subseteq V$ such that $S_1 \ne \emptyset$ and $S_2 \not=
\emptyset$ and $S_1 \cup S_2 = V$. Let there be no vertex $u \in S_2$
such that $N[u] \cap S_1 = S_1$. Let $\opt = \boxi(G)$.  Let $t$ be
a fixed positive integer such that $1 \le t \le |S_1|$.  Let $n_t
= n_t(S_1,S_2)$.  Let $t^* = |S_2| - 2\opt (|S_2| - n_t)$.
Let $m^{*} = m_{t^*}(S_2, S_1)$, for $t^* > 0$. For $t* \le 0$, we define $m^* = \infty$.
  Then
$$
  \boxi(G) \ge \frac{|S_2|}{2 \left( (|S_2|- n_t)   + (t-1) \frac {t^*} {m^*}  \right) } \enspace .
$$
\end{theorem}

\begin{remark}\label{rem:remSecondMethod}
Before moving on to the proof, we would like to bring to
notice certain observations regarding the lower bound.
\begin{enumerate}

\item $m_t>0$ for $1 \le t \le |S_1|$.   This follows from the fact that
there is no $u\in S_2$ such that $N[u]\cap S_1=S_1$.
If $t < 0$, by definition $m^* = \infty$. 
Therefore  $m^*>0$, in all cases.

\item $(|S_2| - n_t)+  (t-1) \frac  {t^*} {m^*}  > 0$. If $t=1$, then,
$(|S_2|-n_t)>0$.  This is because, if $n_1=|S_2|$, it implies that
every vertex in $S_1$ is connected to every vertex in $S_2$, which
contradicts the assumption that there is no $u\in S_2$ such that
$N[u]\cap S_1=S_1$. If $t>1$, then $(t-1)\frac {t^*}{m^*}=0$ only
if $t^* \le 0$. But, when $(|S_2|-n_t)=0$,
$t^*= |S_2| + 2\opt (|S_2| - n_t) = |S_2|>0$.

\ignore { 
\item The maximum value of 
$\frac{|S_2|}{2\left(|S_2|-n_t)+ (t-1) \frac {t^*} {m^*} \right)}$ can be
bounded above by $\max\left(\frac{|S_1|}{2},\frac{|S_2|}{2}\right)$. To see
this, we consider the following two cases:
\begin{enumerate}
\item Suppose $|S_2|-n_t>0$. Then $|S_2|-n_t \ge 1$ and thus the
  expression is at most $\frac {|S_2|} {2}$. 
\item Suppose $|S_2|-n_t=0$. From the above discussion $t>1$ and
$t^*=|S_2|$. Since, $m* \le |S_1|$, 
the expression is at most  $|S_1|/2$.
\end{enumerate}
}

\item We can assume that $t^* > 0$ because if $t^* \le 0$, the lemma easily
holds as follows: From the definition of $t^*$,
we have $|S_2| - 2\opt (|S_2| - n_t)\le0$. Clearly $(|S_2|-n_t)\ne0$
in this case. Therefore,
\[
\opt \ge \frac{|S_2|}{2(|S_2|-n_t)}\ge \frac{|S_2| }{2 \left( ( |S_2|- n_t)   +
(t-1) \frac {t^{*}} {m^*}\right)} 
\]

recalling that   $m^* = \infty$,  for $t^* \le 0$, and therefore 
the second term of the expression in the denominator becomes $0$ .
\end{enumerate}
\end{remark}

\begin{proof}[of Theorem \ref{thm:thm_main}]
Let $I_1, \ldots I_{\opt}$ be $\opt$ interval graphs such that $E(G)
= E(I_1) \cap \cdots \cap E(I_{\opt})$. For each interval graph $I_i$,
let $f_i$ be an interval representation with distinct end-points. We denote
$$
\lt{v}_i=l_{f_i}(v)\mbox{ and }\rt{v}_i = r_{f_i}(v) \enspace ,
$$
Let $x_i \in S_1$ be the vertex such that $ | \{ w \in S_1 ~|~ \lt{w}_i >
\lt{x}_i \} |  =  t -1$.  Similarly, let $y_i\in S_1$ be the vertex such
that $| \{ w \in S_1 ~|~ \rt{w}_i < \rt{y}_i \} | = t -1$. Define
$$
\lt{\mathbf{M}}_i =  \lt{x_i}~\mbox{and}~\rt{\mathbf{m}}_i =  \rt{y_i} \enspace 
$$
and let 
$$
\rt{X}_i = \{ v \in S_2 : \lt{v}_i > \rt{\mathbf{m}}_i\} ~~~~\mbox{and}~~~~ \lt{Y}_i = \{ v \in S_2 : \rt{v}_i < \lt{\mathbf{M}}_i\}  \enspace .
$$

\begin{clm} \label{clm:cl_size}
For all $i \in \{1, \ldots, \opt\}$, $|\rt{X}_i| \le |S_2|  - n_t $ ~and~ $|\lt{Y}_i| \le |S_2| - n_t $.
\end{clm}

\begin{proof}
Let $T_i = \{v \in S_1 ~|~ \rt{v_i} \le \rt{\mathbf{m}}_i \}$. 
From the definition of  $\rt{\mathbf{m}}_i$, we have $|T_i|=t$. 
It follows from Lemma \ref{lem_basic} that for each $v \in N[T_i]$, $\lt{v}_i < \rt{\mathbf{m}}_i$.
Thus, if $u \in N[T_i] \cap S_2$ then $u \notin \rt{X}_i$.
Hence $\rt{X}_i \subseteq S_2 - \left( N[T_i] \cap S_2 \right)$.
Recalling the definition of
$n_t = n_t(S_1,S_2)$, we have $|N[T_i] \cap S_2| \ge
n_t$. 
It follows that  $|\rt{X}_i| \le |S_2| - n_t$.
Using similar arguments, it follows that $|\lt{Y}_i| \le |S_2|   -
n_t$. \qed
\end{proof}
Let $Z = S_2 -\bigcup_{i \in \{1, \ldots, \opt\}} \rt{X}_i\cup\lt{Y}_i$, 
i.e. 
\[
Z = \left\{ u \in S_2 ~|~ \lt{u}_i < \rt{\mathbf{m}}_i \mbox{~and~}
\rt{u}_i > \lt{\mathbf{M}}_i \mbox{~for all~} i \in \{1, \ldots,
\opt\}\right\}.
\]

\begin{clm} \label{clm:cl_z}
 $|Z| \ge t^*>0$.
\end{clm}
\begin{proof}
First, we recall the assumption that $t^*>0$ (see Remark
\ref{rem:remSecondMethod}, point 3).  Using the definition of $Z$ and
Claim \ref{clm:cl_size}, it follows that $|Z| \ge |S_2| - \sum_{i \in \{1,
\ldots, \opt\}} |X_i^+| - \sum_{i \in \{1, \ldots, \opt\}} |Y_i^-| \ge
|S_2| - 2 \opt (|S_2|  - n_t) 
= t^*$.\qed
\end{proof}

\begin{clm} \label{clm:cl_each_inter}
For each interval graph $I_i$ , $i \in \{1, \ldots, \opt\}$, $$| \{v
\in S_1 ~|~ \exists u \in Z \mbox{~and~} uv \notin E(I_i) \} | \le 2 (
t - 1 ).$$
\end{clm}
\begin{proof}
Fix any interval graph $I_i$.  Let $R_i = \{ u \in S_1 ~|~ \rt{u}_i
\ge \rt{\mathbf{m}}_i \}$ and let $L_i = \{ u \in S_1 ~|~ \lt{u}_i \le
\lt{\mathbf{M}}_i \}$.  Recalling the definition of $\rt{\mathbf{m}}_i$
and $\lt{\mathbf{M}}_i$, it follows that  $|R_i| = |L_i| = |S_1| -
(t-1)$. From the definition of $Z$, we know that for each vertex $v \in
Z$, $\lt{v}_i < \rt{\mathbf{m}}_i$ and $\rt{v}_i > \lt{\mathbf{M}}_i$. If
vertex $u \in L_i \cap R_i$ then $\rt{u}_i \ge \rt{\mathbf{m}}_i>\lt{v}_i$ and
$\lt{u}_i \le \lt{\mathbf{M}}_i<\rt{v}_i$.  It follows from Lemma \ref{lem_basic}
that $uv \in E(I_i)$ for all $v\in Z$.  In other words, $ \{v \in S_1 ~|~ \exists u \in
Z \mbox{~and~} uv \notin E(I_i) \}  \subseteq S_1 - ( R_i \cap L_i) =
( S_1 - R_i )  \cup (S_1 - L_i)$.  Since $|S_1 - R_i| = t -1$ and $|S_1 -
L_i| = t -1$, we obtain the claimed bound.\qed
\end{proof}
From the definition of $m_t$ it follows that $|N'(Z,\comp{G})
\cap S_1| \ge m_{|Z|}(S_2, S_1) $. Recalling that $|Z|\ge t^*$
(from Claim \ref{clm:cl_z})  we obtain
\begin{equation} \label{eq_two}
|N'(Z,\comp{G}) \cap S_1| \ge m_{|Z|} (S_2, S_1)  \ge m_{t^*}
(S_2, S_1)  = m^*.
\end{equation}
Clearly $N'(Z,\comp{G}) \cap S_1 = \{v \in S_1 ~|~ \exists u \in Z
\mbox{~and~} uv \notin E(G) \}$.  Since each missing edge in $G$ should be
missing in at least one of the $\opt$ interval graphs, it follows from
Claim \ref{clm:cl_each_inter} that $|N'(Z,\comp{G}) \cap S_1| \le 2 \opt
(t-1)$.  Combining this with inequality (\ref{eq_two}), we obtain
$2 \opt(t-1) \ge |N'(Z,\comp{G}) \cap S_1| \ge m^*.$
Recalling, $t^* = |S_2|  - 2 \opt( |S_2| - n_t)$, we get:
\begin{equation} \label{eq_lower}
2 ( t - 1 )  \opt t^*  \ge  \ m^* (|S_2| - 2\opt (|S_2| - n_t)).
\end{equation}
Rearranging the terms, it follows that
\begin{equation} \label{eq_3}
   2 \opt \left(  m^*  (|S_2|- n_t)  + (t-1) t^*  \right)   \ge
   m^* |S_2|. 
\end{equation}
Recalling that $m^* >0$ and thus  $ |S_2|-n_t+(t-1) \frac {t^*} {m^*} >0$
(see Remark \ref{rem:remSecondMethod}, points 1 and 2), on rearranging
the above inequality we get the required lower bound.  \bqed
\end{proof}

\subsection {Consequences}\label{sec:consMainThm}
In a graph $G$, a vertex is called a \emph{universal vertex} if it is
adjacent to all other vertices.
\begin{theorem}\label{thm:universal}
Let $G$ be a non-complete graph with minimum degree $\delta$.
Let $n_u$ be the number of universal vertices in $G$. Then
$$
	\boxi(G) \ge \frac{n - n_u}{ 2 ( n - \delta - 1) } \enspace .
$$
\end{theorem}
\begin{proof}
Let $U$ be the set of universal vertices in $G$. We have $|U| = n_u$.
Let $G'$ be the graph induced by $V - U$ on $G$.
Let $\delta'$ be the minimum degree of $G'$. It is easy to see
that $\delta' = \delta - n_u$.
We apply Theorem \ref{thm:thm_main} on $G'$ as follows:
take $S_1 = S_2 = V - U$ and fix $t = 1$.
Note that, we can apply Theorem \ref{thm:thm_main} since there are no
universal vertices in $G'$. Clearly $n_t = n_1(S_1, S_2) =  \delta' + 1$,
and thus $|S_2| - n_t = (n-n_u )- (\delta-n_u +1)  = n - \delta - 1$. 
Hence, from Theorem
\ref{thm:thm_main}
we obtain that 
$\boxi(G') \ge    \frac{n - n_u}{ 2(n - \delta -
1)}$.\bqed
\end{proof}
For a $(n-k-1)$-regular graph, this result is same as Theorem
\ref{thm:nk2}. Note however that for general graphs, the above result
is better than the one mentioned in Section \ref{sec:immAppMethod1},
namely $\boxi(G)\ge \frac{n\delta(\comp{G})}{2\Delta(\comp{G})^2}$.

\begin{theorem}\label{thm:biUniversal}
Let $G$ be a non-complete bipartite graph with two parts $A$ and
$B$. Let $u_B$ be the number  of vertices in $B$ that are adjacent to
all vertices in $A$. Let $\delta_A$ be the minimum degree of the vertices
in $A$ . Then,
$$
	\boxi(G) \ge \frac{|B| - u_B }{ 2 ( |B| - \delta_A) } \enspace .
$$
\end{theorem}
\begin{proof}
Let $U_B \subseteq B$ be the set of vertices that are adjacent to
all the vertices in $A$. We have $|U_B| = u_B$. Let $G'$ be the graph
induced by $A \cup (B - U_B)$ on $G$.  Apply Theorem \ref{thm:thm_main}  on
$G'$ as follows: fix $S_1 = A$, $S_2 = B - U_B$ and $t = 1$. Note that
Theorem \ref{thm:thm_main} can be applied since there is no vertex $u\in
B-U_B$ that is adjacent to all vertices in $A$. Clearly,
the minimum degree of the vertices in $S_1$ is $\delta'_A = \delta_A -
u_B$.  We have $n_t = n_1(S_1, S_2) = \delta'_A$,
and thus $|S_2| - n_t = |B| - \delta_A$. Applying Theorem
\ref{thm:thm_main} we obtain that $\boxi(G') \ge  
\frac{|B| - u_B}{2(|B| - \delta_A)}$.\bqed
\end{proof}

Note that in the case of  regular balanced bipartite graphs,
the above theorem is same as Theorem  \ref{thm:bipnk2}, but
in the non-regular case the above theorem gives a better bound
than what can be obtained using the first method.

\begin{theorem}\label{thm:tExp}
Let $G$ be a graph on $n$ vertices with minimum degree
$\delta$ and maximum degree $\Delta<n-1$. If $\overline G$ 
does not have  $K_{t,t}$ as a subgraph,
 then
$$
\boxi(G) \ge \frac{n (n - \Delta-1)}{2(t-1)\left[(n - \Delta-1) + (n -
\delta-1)\right]} \enspace .
$$
\end{theorem}
\begin{proof}
We apply Theorem \ref{thm:thm_main} on $G$ as follows. Fix $S_1 = S_2
= V$.  It is easy to see that if $\overline {G}$ does not have
  $K_{t,t}$ as a subgraph, then $n -
n_t  \le t-1$. 
Clearly we have   $t^* (n-\Delta-1) \le m^* (n-\delta -1)$.
Substituting for $n-n_t$ and $ \frac {t^*} {m^*}$
in Theorem \ref{thm:thm_main} we obtain the result. \bqed
\end{proof}
In Theorem \ref{thm:nkExamples} of Section \ref{sec:immAppMethod1}, we
obtained lower bounds for the boxicity of regular co-planar graphs and
regular $C_4$-free graphs based on the absence of $K_{3,3}$ and $K_{2,2}$
respectively as subgraphs in their complements. We observe that these
graphs are specific subclasses of regular graphs whose complements
does not have $K_{t,t}$ as a subgraph. The same results
can be obtained using the following corollary to Theorem \ref{thm:tExp}.

\begin{corollary}
Let $G$ be a regular non-complete graph such that $\overline {G}$ does
not have  
$K_{t,t}$ as a subgraph. Then 
$$
\boxi(G) \ge \frac{n}{4(t-1)} \enspace .
$$
\end{corollary}

\section*{Acknowledgments}
We thank Joel Friedman for personal communication regarding his proof of
Alon's conjecture. We also thank the anonymous referee for the comments
which improved the exposition.

\end{document}